\documentclass[11pt]{article}

\usepackage{amssymb,amsmath,accents,amsfonts,amsthm}

\usepackage{mathtools}
\mathtoolsset{showonlyrefs}

\usepackage[font=small,labelfont=bf]{caption}
\usepackage{makeidx}
\makeindex

\usepackage{xspace}

\usepackage[verbose,letterpaper,tmargin=2.54cm,bmargin=2.54cm,lmargin=2.54cm,rmargin=2.54cm]{geometry}

\usepackage{graphicx} 
\usepackage{empheq} 
\usepackage{color,framed}
\usepackage{float}
\usepackage{epsfig} 
\usepackage{pdfpages}

\usepackage[utf8x]{inputenc}
\usepackage[english]{babel}
\usepackage[T1]{fontenc}

\usepackage[bookmarks,bookmarksnumbered,colorlinks=true,urlcolor=black,linkcolor=black,citecolor=blue]{hyperref}
\usepackage{memhfixc}

\usepackage[numbers,comma,sort&compress]{natbib}

\usepackage[noblocks]{authblk}

\newcommand*{\bL}{{\bf L}}

\newcommand*{\bS}{{\bf S}}
\newcommand*{\bT}{{\bf T}}

\newcommand*{\bW}{{\bf W}}

\newcommand*{\ba}{{\bf a}}

\newcommand*{\be}{{\bf e}}
\newcommand*{\bff}{{\bf f}}

\newcommand*{\bx}{{\bf x}}
\newcommand*{\by}{{\bf y}}
\newcommand*{\bz}{{\bf z}}

\newcommand*{\rd}{\mathrm{d}}

\newcommand*{\rp}{\mathrm{p}}

\newcommand*{\rT}{\mathrm{T}}

\newcommand*{\A}{\mathbb{A}}

\newcommand*{\D}{\mathbb{D}}
\newcommand*{\E}{\mathbb{E}}

\renewcommand*{\H}{\mathbb{H}}

\newcommand*{\J}{\mathbb{J}}

\renewcommand*{\P}{\mathbb{P}}

\newcommand*{\R}{\mathbb{R}}
\renewcommand*{\SS}{\mathbb{S}}
\newcommand*{\T}{\mathbb{T}}
\newcommand*{\U}{\mathbb{U}}
\newcommand*{\V}{\mathbb{V}}

\newcommand*{\Lc}{\mathcal{L}}

\newcommand*{\Oc}{\mathcal{O}}

\newcommand*{\kvec}[2][]{#1\left(\begin{array}{c}#2#1\end{array}\right)}

\newcommand*{\abs}[2][]{#1\left\vert#2#1\right\vert}

\newcommand*{\norm}[2][]{#1\left\Vert#2#1\right\Vert}

\newcommand*{\parenth}[2][]{#1\left(#2#1\right)}

\newcommand*{\curly}[2][]{#1\left\{#2#1\right\}}

\newcommand*{\smcurly}[2][]{#1\big\{#2#1\big\}}

\newcommand*{\scalar}[2][]{#1\left\langle#2#1\right\rangle}

\newcommand*{\smscalar}[2][]{#1\big\langle#2#1\big\rangle}

\newcommand*{\block}[2][]{#1\left[#2#1\right]}

\newcommand*{\smblock}[2][]{#1\big[#2#1\big]}
\newcommand*{\mblock}[2][]{#1\Big[#2#1\Big]}

\newcommand*{\parmatrix}[2][]{\begin{pmatrix}#2\end{pmatrix}}

\newcommand*{\limwith}[2][]{#1\underset{#2\to\infty}{\lim}#1}
\DeclareMathOperator{\vspan}{span}
\DeclareMathOperator{\Id}{Id}

\DeclareMathOperator{\ACV}{ACV}

\DeclareMathOperator{\SO}{SO}
\DeclareMathOperator{\PS}{PS}
\DeclareMathOperator{\NSR}{NSR}

\newcommand{\yr}{\rm yr}
\newcommand*{\bdot}[1]{\accentset{\bullet}{#1}}

\newcommand*{\lequation}[2]{\begin{equation}\label{#1}\begin{split}#2\end{split}\end{equation}} 
\newcommand*{\lfigure}[3]{\begin{figure}[h]\centering#3\caption{\small #2}\label{#1}\end{figure}} 

\title{Stable limit cycles perturbed by noise}
\author[1]{\small Stilianos Louca}
\affil[1]{\emph{Institute of Applied Mathematics, University of British Columbia, 121-1984 Mathematics Road, Vancouver, BC, V6T1Z2 Canada, \href{mailto:louca@math.ubc.ca}{louca@math.ubc.ca}}}

\begin{document}

\maketitle

\begin{abstract}
Many physical and biological systems exhibit intrinsic cyclic dynamics that are altered by random external perturbations. We examine continuous-time autonomous dynamical systems exhibiting a stable limit cycle, perturbed by additive Gaussian white noise. We derive a formal approximation for the dynamics of sample paths that stay close to the limit cycle, in terms of a \emph{phase} coordinate and a \emph{deviation} perpendicular to the limit cycle. 
To leading order in the deviation, the phase advances at the deterministic speed superimposed by a Brownian-motion-like drift. 
The deviation itself takes the form of an $(n-1)$-dimensional Ornstein-Uhlenbeck process. 
We apply these results to the case of limit cycles emerging through a supercritical Hopf bifurcation, which is widespread in ecological and epidemiological models. We derive approximation formulas for the system's stationary autocovariance and power spectral density. The latter two reflect the effects of perturbations on the temporal coherence and spectral bandwidth of perturbed limit cycles. We verify our results using numerical simulations and exemplify their application to the El Ni\~no Southern Oscillation.
\end{abstract}


\section{Introduction}
Many deterministic mathematical models exhibit stable limit cycles that are used to describe a wide spectrum of natural phenomena, ranging from population cycles \citep{May1972Limit,Scheffer1991Fish,Nisbet2004Modelling}, chemical oscillations \citep{Schnakenberg1979Simple}, geophysical cycles \citep{Payne1995Limit} to periodic epidemic outbreaks \citep{Hethcote1981Nonlinear,Alexander2004Periodicity} and genetic oscillators \citep{Vilar2002Mechanisms}. Systems exhibiting intrinsically emerging oscillations are typically subject to external perturbations (e.g. the random change of an environmental parameter) or internal stochasticity (e.g. due to finite population sizes).
Models are often extended to include such random factors by adding a stochastic term to the deterministic equations, typically in the form of additive Gaussian white noise \citep{Rozenfeld2001On-the-influence,Wang2005Internal,Ji2009Analysis}. 

The effects of noise on limit cycles can range from frequency shifts \citep{Gang1993Stochastic} to an increased spectral bandwidth and a rapid decay of the cycle's autocorrelation \citep{Burgers1999The-El-Nino,Nisbet2004Modelling}.
Quantifying these properties of \emph{noisy limit cycles} is unavoidable if we want to compare the predictions of our models to recorded time series. Such a comparison is complicated by the fact that noise can also induce decoherent oscillations, in cases where the deterministic model predicts a damped oscillation towards a stable fixed point \citep{Aparicio2001Sustained,Allen2001Quasi-cycles,McKane2005Predator-Prey,McKane2007Amplified,Pineda-Krch2007A-tale}. 
\citet{Baxendale2011Sustained} show that under certain approximations, such \emph{quasi-cycles} are described by circular orbits whose radius behaves like an Ornstein-Uhlenbeck process with zero mean \citep{Uhlenbeck1930On-the-Theory,Gillespie1992Markov}. \citet{Thompson2012Stochastic} extend these results to pairs of quasi-cycles sustained by a common noise source.
 \citet{Boland2008How-limit} investigate limit cycles of the two-dimensional Brusselator system subject to noise, by describing the dynamics along the directions tangential and normal to the limit cycle. On a more abstract setting, \citet{Teramae2009Stochastic} and \citet{Goldobin2010Dynamics} derive  phase equations for noisy limit cycles.

In the present article, we consider arbitrary autonomous dynamical systems described by an ordinary differential equation exhibiting a stable limit cycle, perturbed by additive Gaussian white noise. We formally derive an approximation for the dynamics of sample paths that stay close to the limit cycle. We formulate these dynamics as a set of coupled stochastic differential equations (SDE) for the \emph{phase} (or \emph{longitudinal}) coordinate along the limit cycle and a set of complementary coordinates collectively referred to as \emph{deviation}. A similar decomposition was used by \citet{Kurrer1991Effect} to investigate the effects of noise on the Bonhoeffer-van der Pol nonlinear oscillator.
In analogy to previous work \citep{Teramae2009Stochastic}, we find that the phase advances at the deterministic speed superimposed by a Brownian-motion-like drift. 
The deviation takes the form of an $(n-1)$-dimensional Ornstein-Uhlenbeck process \citep{Uhlenbeck1930On-the-Theory,Gillespie1992Markov}, extending the work of \citet{Baxendale2011Sustained} on quasi-cycles. We apply these results to systems in the normal form of the supercritical Hopf bifurcation. The latter describes the emergence of a limit cycle around an unstable focus and appears in numerous ecological and epidemiological models \citep{Rosenzweig1963Graphical,Greenhalgh1997Hopf,Fussmann2000Crossing,Pujo-Menjouet2004Contribution}.
We give explicit approximation formulas for the long-term autocovariance and power spectral density of the components of such a process. These formulas reproduce the widely observed frequency shift and decoherence of noisy limit cycles, and allow a direct comparison of models exhibiting cyclic dynamics to real time series. We asses the fidelity of the derived approximations using numerical simulations and exemplify their use for the El Ni\~no Southern Oscillation, a well known but poorly understood cyclic climatic phenomenon \citep{Trenberth1997The-Definition}.

\section{Linear approximation of noisy limit cycles}
\label{SECTION:LIN_APPROX_NOISY_LC}
Our starting point is a smooth $n$-dimensional dynamical system $\frac{d\by}{dt}=\bff(\by)$, exhibiting a stable limit cycle $\Lc$. With additive Gaussian white noise, the dynamics take the form of an It\^o SDE \citep{Oksendal2003Stochastic}
\lequation{EQ:LIN_APPROX_NOISY_LC01}{
	d\by=\bff(\by)\ dt + \SS\ d\bW,
}
where $\bW$ shall be an $n$-dimensional Wiener process (or \emph{standard Brownian motion}) with uncorrelated components and $\SS\in\R^{n\times n}$ is some matrix.

Let any sample path $\by(t)$ of the SDE \eqref{EQ:LIN_APPROX_NOISY_LC01} close to $\Lc$ be split into a \emph{longitudinal} and a \emph{lateral} component, $\by(t) = \by_{\rm lon}(t) + \by_{\rm lat}(t)$, where $\by_{\rm lon}(t)\in \Lc$ is a point on the limit cycle and $\by_{\rm lat}(t)$ is perpendicular to the local tangent. 
Let $\bL(t)$ be a parameterization of the deterministic limit cycle, i.e. such that $\bL(0) = \by_{\rm lon}(0)$ and $\frac{d}{dt}\bL = \bff(\bL)$. Let $\bT(t)$ be the normalized tangent at $\bL(t)$ in the direction of motion and let $P(t)$ be the hyperplane perpendicular to $\bT(t)$. Abbreviate $\bT_0=\bT(0)$ and $P_0=P(0)$. The pair $(\bT(t),P(t))$ defines a comoving frame along the limit cycle, similar to the Frenet frame in 3 dimensions \citep{Gibson2001Elementary}. Let $\U(t)\in\SO(n)$ be an orthogonal matrix, depending smoothly on time $t$, such that $\bT(t)=\U(t)\bT_0$, $P(t)=\U(t)P_0$ and $\frac{d}{dt}\U(t)\bz\perp P(t)$ for all $\bz\in P_0$. We refer to appendix \ref{LEMMA:EX_CERT_ORTH_TRAFO} for a constructive proof of existence. $\U(0)$ can be chosen to be the identity matrix $\Id$, in which case $\U(t)$ would be unique. For example, in two dimensions $\U(t)$ would be the rotation that maps $\bT_0$ to $\bT(t)$. See figure \ref{FIG:PROOFCOORDINATES01} for an illustration.
\lfigure{FIG:PROOFCOORDINATES01}{Representation of the state $\by$ by a phase $\tau$ and a deviation $\bz\in P_0$, illustrated for a 3 dimensional phase space. The phase $\tau$ advances such that the plane perpendicular to the tangent $\bT(\tau)$, $P(\tau)$, includes $\by$. The deviation $\bz$ is defined in the fixed plane $P_0$.}{
\includegraphics[width=20em]{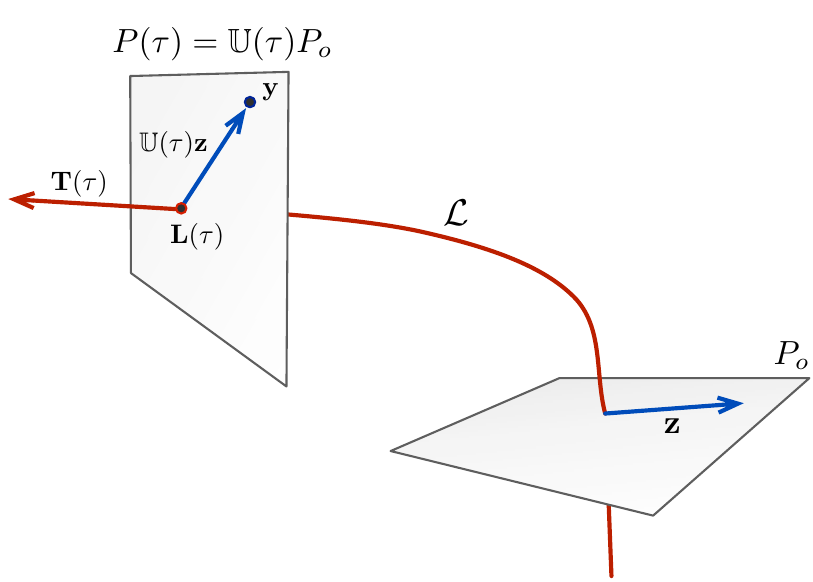}
}
Denote $\V(t)=\frac{d}{dt}\U(t)$. Let $\P(t)$ and $\T(t)$ be the orthogonal projections onto $P(t)$ and the linear span of $\bT(t)$, respectively. Note that $\bT(t),\T(t),P(t),\P(t),\U(t)$ and $\V(t)$ are solely determined by the kinematics of the deterministic trajectory $\bL$. Let $\J(t)=(\nabla\bff)|_{\bL(t)}$ be the Jacobian of the deterministic dynamics around $\bL(t)$.

We make the ansatz $\by_{\rm lon}(t) = \bL(\tau(t))$ and $\by_{\rm lat}(t) = \U(\tau(t))\bz(t)$ for suitable $\tau(t)\in\R$ and $\bz(t)\in P_0$. Upon choice of a basis in the $(n-1)$-dimensional hyperplane $P_0$, $\bz(t)$ can be described by $(n-1)$ independent variables. Hence, $\tau(t)$ and $\bz(t)$ are to be seen as a set of new coordinates for $\by(t)$ along and perpendicular to the limit cycle, which we shall refer to as \emph{phase} and \emph{deviation}, respectively. In these coordinates the kinetics formally take the form
\lequation{EQ:LIN_APPROX_NOISY_LC04}{
	\frac{d\by}{dt} = \frac{d\tau}{dt}\bff(\bL(\tau)) + \frac{d\tau}{dt}\V(\tau)\bz + \U(\tau)\frac{d\bz}{dt},
}
while we omit to show the explicit dependence on $t$ for brevity. On the other hand, linearising the dynamics in the proximity of the limit cycle yields the formal approximation
\lequation{EQ:LIN_APPROX_NOISY_LC06}{
	\frac{d\by}{dt} = \bff(\by_{\rm lon}) + \J(\tau)\by_{\rm lat} + \SS \ \frac{d\bW}{dt} + \Oc(H\norm{\by_{\rm lat}}^2),
}
where $H$ is the bound of the second derivative of $\bff$ (when considered as a bilinear operator $\R^n\times \R^n\to\R^n$), maximized along the entire deterministic limit cycle.
Combining \eqref{EQ:LIN_APPROX_NOISY_LC04} with \eqref{EQ:LIN_APPROX_NOISY_LC06} gives
\lequation{EQ:LIN_APPROX_NOISY_LC08}{
	&\parenth{\frac{d\tau}{dt}-1}\block{\bff(\bL(\tau)) + \V(\tau)\bz}\\
	=& \J(\tau)\U(\tau)\bz + \SS \ \frac{d\bW}{dt} - \U(\tau)\frac{d\bz}{dt}- \V(\tau)\bz + \Oc(H\norm{\bz}^2).
}
Equation \eqref{EQ:LIN_APPROX_NOISY_LC08} can be split into the lateral part
\lequation{EQ:LIN_APPROX_NOISY_LC10}{
	\frac{d\bz}{dt} = \U^\rT(\tau)\P(\tau)\block{\J(\tau)\U(\tau)\bz + \SS \ \frac{d\bW}{dt}} + \Oc(H\norm{\bz}^2)
}
and the longitudinal part
\lequation{EQ:LIN_APPROX_NOISY_LC12}{
	&\parenth{\frac{d\tau}{dt}-1}\block{\bff(\bL(\tau)) + \V(\tau)\bz} \\
	=& \T(\tau)\block{\J(\tau)\U(\tau)\bz + \SS \ \frac{d\bW}{dt}} - \V(\tau)\bz + \Oc(H\norm{\bz}^2),
}
while we used the fact that $\V(t)\bz\perp P(t)$. Equation \eqref{EQ:LIN_APPROX_NOISY_LC12} can be written as
\lequation{EQ:LIN_APPROX_NOISY_LC14}{
	\frac{d\tau}{dt} =& 1+ \frac{\scalar{\J(\tau)\U(\tau)\bz - \V(\tau)\bz + \SS \ \tfrac{d\bW}{dt}, \bT(\tau)}}{\smscalar{\bff(\bL(\tau)) + \V(\tau)\bz,\bT(\tau)}} \\
	&+ \Oc\parenth{\frac{H\norm{\bz}^2}{\norm{\bff(\bL(\tau))} - \norm{\V(\tau)\bz}}}.
}
Note that
\lequation{EQ:LIN_APPROX_NOISY_LC14B}{
	\norm{\V(\tau)\bz} = \norm{\bff(\bL(\tau))}\cdot \norm{\frac{d\U}{ds}(\tau)\bz},
}
where $s(t)$ is the arc-length parameterization along the limit cycle. Furthermore, $\norm{\frac{d\U}{ds}(\tau)}\leq \norm{\frac{d\bT}{ds}(\tau)}$ (see \ref{APPENDIX:NORM_ROTATIONS} for a proof). Hence, for small $\bz$, \eqref{EQ:LIN_APPROX_NOISY_LC14} can be written as
\lequation{EQ:LIN_APPROX_NOISY_LC14D}{
	\frac{d\tau}{dt} =& 1+ \frac{\scalar{\J(\tau)\U(\tau)\bz - \V(\tau)\bz + \SS \ \tfrac{d\bW}{dt}, \bT(\tau)}}{\smscalar{\bff(\bL(\tau)) + \V(\tau)\bz,\bT(\tau)}} \\
	&+ \Oc\parenth{\frac{H}{\norm{\bff(\bL(\tau))}}\cdot\frac{\norm{\bz}^2}{1- \kappa(\tau)\norm{\bz}}},
}
where $\kappa=\norm{\frac{d\bT}{ds}}$ is the curvature of the deterministic limit cycle.
Summarizing, we can decompose sample paths of \eqref{EQ:LIN_APPROX_NOISY_LC01} as
\lequation{EQ:LIN_APPROX_NOISY_LC15}{
	\by= \bL(\tau) + \U(\tau)\bz,
}
where $\bz$ and $\tau$ are correlated stochastic processes whose sample paths approximately solve \eqref{EQ:LIN_APPROX_NOISY_LC10} and \eqref{EQ:LIN_APPROX_NOISY_LC14D}, respectively.

Due to the local stability of the limit cycle, $\U^\rT(t)\P(t)\J(t) \U(t)$ is a stable linear operator on $P_0$. Therefore, \eqref{EQ:LIN_APPROX_NOISY_LC10} describes an Ornstein-Uhlenbeck process on the hyperplane $P_0$, with a possibly stochastic noise tensor and Jacobian \citep{Gardiner1985Handbook}. Intuitively, deviations of sample paths from the limit cycle are the result of fluctuations acting against local stabilising dynamics \citep{Ali1999On-the-local}, characterized predominantly by the limit cycle's Lyapunov exponent. This is in accordance with results by \citet{DeVille2011Stability}, who showed that for weak noise limit cycles in Hopf normal form have a negative Lyapunov exponent.
On the other hand, to leading order in $\bz$, the phase $\tau$ advances at the deterministic rate modulated by additive white noise, therefore exhibiting a Brownian-motion-like drift away from its deterministic value. As such, \eqref{EQ:LIN_APPROX_NOISY_LC14D} differs fundamentally from \eqref{EQ:LIN_APPROX_NOISY_LC10}, since noise-induced changes of the phase are not reversed by the deterministic dynamics.

Without loss of generality let $P_0=\curly{(x_1,..,x_n)\in\R^n : x_n=0}$, so that $\bz=(\bz_0,0)$ for some $\bz_0\in\R^{n-1}$, and let $\Pi_0:\R^n\to\R^{n-1}$ denote the projection to the first $(n-1)$ components. Suppose that the noise $\SS \ d\bW$ is isotropic with uncorrelated components, i.e. $\SS^\rT\SS=\sigma^2\cdot\Id$ for some scalar $\sigma\in\R$. Then \eqref{EQ:LIN_APPROX_NOISY_LC10} and \eqref{EQ:LIN_APPROX_NOISY_LC14D} can be written to leading order as
\lequation{EQ:LIN_APPROX_NOISY_LC16}{
	&d\bz_0 \approx \J_0(\tau)\bz_0\ dt + \sigma\ d\bW_\rd,\quad d\tau \approx dt + \frac{\sigma\ dW_\rp}{\norm{\bff(\bL(\tau))}},
}
where $\J_0$ is defined by
\lequation{EQ:LIN_APPROX_NOISY_LC16B}{
	\J_0(\tau)\bz_0=\Pi_0\U^\rT(\tau)\P(\tau)\J(\tau) \U(\tau)\parmatrix{\bz_0\\0},
}
and $\bW_\rd$ and $W_\rp$ are Wiener processes of dimension $n-1$ and $1$, respectively. In fact, $\bW_\rd$ and $W_\rp$ are uncorrelated, because the two orthogonal projections $\T(\tau)$ and $\P(\tau)$ split the noise $\SS\ d\bW$ into two uncorrelated processes.

\section{Limit cycles emerging through Hopf bifurcations}
\label{SECTION:LC_HOPF}
In this section, we exemplify the linear approximation \eqref{EQ:LIN_APPROX_NOISY_LC15} for the case where the deterministic limit cycle emerges through a supercritical Hopf bifurcation \citep{Hale1991Dynamics,Strogatz1994Nonlinear}. The ubiquity of studied dynamical systems exhibiting a Hopf bifurcation makes this an ideal illustrative example. 

\subsection{Approximation of sample paths}
\label{SECTION:LC_HOPF_APPROX_SAMPLE_PATHS}
For simplicity, we consider the normal form \cite[\S 3.4]{Guckenheimer1985Nonlinear}
\lequation{EQ:LIN_APPROX_NOISY_LC_HOPF01}{
	d\kvec{x\\ y} =& \parmatrix{\lambda/2 & & -\alpha_0\\\alpha_0 & & \lambda/2}\cdot\kvec{x\\y}\ dt\\
	& + \frac{(x^2+y^2)}{r^2} \kvec{-\lambda x/2 - (\alpha-\alpha_0) y\\-\lambda y/2 + (\alpha-\alpha_0) x}\ dt\\
	& + \sigma \ d\bW,\noeqref{EQ:LIN_APPROX_NOISY_LC_HOPF01}
}
where $\alpha>0$, $r>0$ and $-\lambda<0$ are the limit cycle's angular frequency, radius and Lyapunov exponent, respectively, $\alpha_0>0$ is the system's angular frequency in the proximity of the focus and $\sigma\neq0$. The first row in \eqref{EQ:LIN_APPROX_NOISY_LC_HOPF01} describes the linear dynamics in the proximity of the focus, which has Lyapunov exponent $\lambda/2$. The second row describes the nonlinearities giving rise to the limit cycle. The Wiener process $\bW$ appearing in the third row is assumed to have uncorrelated components. In polar coordinates $(\rho,\varphi)$ the deterministic part of \eqref{EQ:LIN_APPROX_NOISY_LC_HOPF01} reads
\lequation{EQ:LIN_APPROX_NOISY_LC_HOPF01B}{
	\frac{d\rho}{dt} = \frac{\lambda}{2}\rho - \frac{\lambda\rho^3}{2r^2},\quad
	\frac{d\varphi}{dt} = \alpha_0 + (\alpha-\alpha_0)\frac{\rho^2}{r^2}.\noeqref{EQ:LIN_APPROX_NOISY_LC_HOPF01B}
}
The SDE \eqref{EQ:LIN_APPROX_NOISY_LC_HOPF01} describes a stationary process with zero mean.
Modulo arbitrary phase shifts, the limit cycle solution to the deterministic part of \eqref{EQ:LIN_APPROX_NOISY_LC_HOPF01} is given by $\bL(t) = \parenth{r\cos\alpha t, r\sin\alpha t}^\rT$. Calculating $\bT(t),\P(t),\J(t),\U(t)$ and $\V(t)$ is straightforward and one eventually obtains from \eqref{EQ:LIN_APPROX_NOISY_LC10}, \eqref{EQ:LIN_APPROX_NOISY_LC14D} and \eqref{EQ:LIN_APPROX_NOISY_LC15} the approximation
\lequation{EQ:LIN_APPROX_NOISY_LC_HOPF03}{
	\kvec{x(t)\\y(t)} \approx \kvec{(r+z(t))\cdot \cos\alpha \tau(t)\\(r+z(t))\cdot \sin\alpha \tau(t)}
}
for the sample paths, where
\lequation{EQ:LIN_APPROX_NOISY_LC_HOPF05}{
	&dz = -\lambda z\ dt + \sigma\ dW_\rd,\quad 
	d\tau \approx dt + \frac{2z(\alpha-\alpha_0)}{\alpha(r+z)}\ dt + \frac{\sigma\ dW_\rp}{\alpha(r+z)}.\noeqref{EQ:LIN_APPROX_NOISY_LC_HOPF05}
}
The SDE for $z$ describes a classical Ornstein-Uhlenbeck process, suggesting that the stationary expected squared distance from the deterministic limit cycle will be approximately $\sigma^2/(2\lambda)$. In the following, the ratio $\NSR=\sqrt{\sigma^2/(2\lambda)}/r$ shall be referred to as \emph{noise-to-signal ratio} (see section \ref{SECTION:LEAD_ORD_APPR_ACV} for further justification of this name). The 2nd term in the SDE for $\tau$ vanishes if the focal frequency is similar to the limit cycle frequency ($\alpha_0\approx\alpha$), which is a reasonable approximation if the system is close to the Hopf bifurcation point. However, for large limit cycle radii this term is expected to become non-negligible (see numerical tests in section \ref{SECTION:NUMTESTS}).

\subsection{Leading order approximation of the autocovariance}
\label{SECTION:LEAD_ORD_APPR_ACV}
We use \eqref{EQ:LIN_APPROX_NOISY_LC_HOPF03} to derive an approximation for the stationary autocovariance,
\lequation{EQ:LEAD_ORD_APPR_ACV01}{
	\ACV[x](u)=\limwith{t}\E\curly{x(t)x(t+u)},
}
of the component $x$. We focus on single components because we wish to draw an analogy to ecological times series, which are often only available for a few system variables. We consider the SDE \eqref{EQ:LIN_APPROX_NOISY_LC_HOPF05} for $\tau$ up to leading order in $z$. More precisely, we approximate
\lequation{EQ:LC_SUPCRIT_HOPFBIF20}{
	x(t) \approx (r + z)\cdot \cos\varphi,\noeqref{EQ:LC_SUPCRIT_HOPFBIF20}
}
where $\varphi(t) = \alpha t + (\sigma/r)W_\rp$ is a Brownian motion with deterministic drift and $dz = - \lambda z + \sigma dW_\rd$ describes a one-dimensional Ornstein-Uhlenbeck process, independent of $\varphi$. The process $\cos\varphi$ is sometimes referred to as \emph{randomised harmonic process} \cite[\S 1.2.1]{Zhu2013Bounded}. Inserting \eqref{EQ:LC_SUPCRIT_HOPFBIF20} into \eqref{EQ:LEAD_ORD_APPR_ACV01} leads to
\lequation{EQ:LEAD_ORD_APPR_ACV03}{
	\ACV[x](u) =&\limwith{t}\frac{r^2}{2}\E\smcurly{\cos\varphi(t)\cos\varphi(t+u)}\\
	&+\limwith{t}\frac{1}{2}\E\curly{z(t)z(t+s)}\E\smcurly{\cos\varphi(t)\cos\varphi(t+u)}\\
	=&\frac{1}{2}\smblock{r^2 + \ACV[z](u)}\cdot \ACV[\cos\varphi](u),
}
where $\ACV[z](u)=\frac{\sigma^2}{2\lambda}e^{-\lambda\abs{u}}$ is the stationary autocovariance of $z$ \cite[\S 3.3.C]{Gillespie1992Markov}.
The autocovariance $\ACV[\cos\varphi](u)$ is well known \cite[eq. (1.9)]{Zhu2013Bounded}, yielding
\lequation{EQ:LEAD_ORD_APPR_ACV05}{
	\ACV[x](u) = \frac{r^2}{2}\mblock{1 + \NSR^2e^{-\lambda \abs{u}}}\cdot \cos(\alpha u)\cdot e^{-\abs{u}(\sigma/r)^2/2}.\noeqref{EQ:LEAD_ORD_APPR_ACV05}
}
Here, $\NSR=\sqrt{(\sigma/r)^2/(2\lambda)}$ is the noise-to-signal ratio introduced in section \ref{SECTION:LC_HOPF_APPROX_SAMPLE_PATHS}. $\NSR^2$ thus relates the signal variance originating in lateral fluctuations around the limit cycle to the variance caused by the limit cycle itself. $1/\NSR^2$ is comparable to the \emph{signal-to-noise ratio} known from signal processing theory \citep{Johnson2006Signal-to-noise}. The right-most exponentially decaying factor corresponds to a temporal decoherence of the limit cycle \citep{Marathay1982Elements}, due to noise-induced phase drift. This deviation from true periodicity in the presence of \emph{phase noise}, is known as \emph{jitter} in electronic signal theory \citep{Demir2000Phase}.

\subsection{Leading order approximation of the power spectral density}
\label{SECTION:LEAD_ORD_APPR_POW_SPEC}
Similarly to section \ref{SECTION:LEAD_ORD_APPR_ACV}, we use the leading order approximation \eqref{EQ:LC_SUPCRIT_HOPFBIF20} to predict the power spectral density,
\lequation{EQ:LC_SUPCRIT_HOPFBIF21}{
	\PS[x](\omega) = \limwith{T} \E\abs{\frac{1}{\sqrt{T}}\int_0^Tx(t)e^{-i\omega t}\ dt}^2,
}
of $x$. Inserting \eqref{EQ:LC_SUPCRIT_HOPFBIF20} into \eqref{EQ:LC_SUPCRIT_HOPFBIF21} yields
\lequation{EQ:LC_SUPCRIT_HOPFBIF22}{
	\PS[x](\omega) = A^2\PS[\cos\varphi](\omega) + \PS[z\cos\varphi](\omega).
}
The power spectral density of $\cos\varphi$ is well known \cite[\S 5]{Xie2006Dynamic} and given by
\lequation{EQ:LC_SUPCRIT_HOPFBIF23}{
	\PS[\cos\varphi](\omega) =& \frac{1}{2}\int_\R \cos(\alpha u)e^{-(\sigma/r)^2\abs{u}/2}\ du\\
	=& \frac{2(\sigma/r)^2\block{4(\alpha^2+\omega^2) + (\sigma/r)^4}}{\block{4(\alpha-\omega)^2 + (\sigma/r)^4}\block{4(\alpha+\omega)^2 + (\sigma/r)^4}}.
}
It is straightforward to see that similarly,
\lequation{EQ:LC_SUPCRIT_HOPFBIF24}{
	\PS[z\cos\varphi](\omega) = \frac{1}{2}\int_\R \ACV[z](u)\cdot \cos(\alpha u)e^{-(\sigma/r)^2\abs{u}/2}\ du.
}
Evaluating \eqref{EQ:LC_SUPCRIT_HOPFBIF24} yields, together with \eqref{EQ:LC_SUPCRIT_HOPFBIF22} and \eqref{EQ:LC_SUPCRIT_HOPFBIF23}, the power spectral density
\lequation{EQ:LC_SUPCRIT_HOPFBIF27}{
	\PS[x](\omega) =& \frac{2r^2(\sigma/r)^2\block{4(\alpha^2+\omega^2) + (\sigma/r)^4}}{\block{4(\alpha-\omega)^2 +  (\sigma/r)^4}\block{4(\alpha+\omega)^2 + (\sigma/r)^4}}\\
	&+\NSR^2\cdot\frac{2r^2\parenth{(\sigma/r)^2 + 2\lambda}\block{4(\alpha^2+\omega^2) + \parenth{(\sigma/r)^2 + 2\lambda}^2}}{\block{4(\alpha-\omega)^2 + \parenth{(\sigma/r)^2 + 2\lambda}^2}\block{4(\alpha+\omega)^2 + \parenth{(\sigma/r)^2 + 2\lambda}^2}}.\noeqref{EQ:LC_SUPCRIT_HOPFBIF27}
}
The power spectrum \eqref{EQ:LC_SUPCRIT_HOPFBIF27} can also be obtained from the autocovariance \eqref{EQ:LEAD_ORD_APPR_ACV05} using the Wiener-Khintchine theorem \citep{Pollock1999Handbook}.
As suggested by \eqref{EQ:LC_SUPCRIT_HOPFBIF27}, the presence of noise leads to a shift of the spectral peak to a higher frequency than the limit cycle's deterministic frequency. Moreover, even for very stable limit cycles ($\lambda\gg(\sigma/r)^2$, or $\NSR\approx0$), the power spectrum retains a non-vanishing bandwidth, leading to the temporal decoherence expressed by the decaying autocovariance \eqref{EQ:LEAD_ORD_APPR_ACV05} \cite[\S 2.3.1]{Glindemann2011Principles}.

\section{Numerical validation}
\label{SECTION:NUMTESTS}
To test the fidelity of our results, we performed numerical simulations of the exact system \eqref{EQ:LIN_APPROX_NOISY_LC_HOPF01} and its linear approximation \eqref{EQ:LIN_APPROX_NOISY_LC_HOPF03} over a wide parameter range. We refer to \ref{APPENDIX:DETAILS_NUM_METHODS} for technical details. The linear approximation is found to have sample paths and a probability distribution that are similar to the exact system, provided that noise is sufficiently weak ($\NSR\lesssim 0.1$) (figures \ref{FIG:XYHOPF_TRAJ01}a,b,c). For stronger noise ($\NSR\gtrsim 0.5$), the linear approximation has a distribution with heavier tails than expected, as well as a stronger peak in the region enclosed by the limit cycle (figures \ref{FIG:XYHOPF_TRAJ01}d,e,f). 
The presence of strong outliers in the linear approximation is due to the linear term $-\lambda z$ in (\ref{EQ:LIN_APPROX_NOISY_LC_HOPF05}), which underestimates the force of attraction by the limit cycle in the outer region, compared to the cubic term $\propto -\rho^3$ in the original system (\ref{EQ:LIN_APPROX_NOISY_LC_HOPF01B}). The false peak in the inner region is due to the fact that in the linear approximation, trajectories with inward deviations exceeding the cycle radius ($z<-r$) have to pass through the origin ($z=-r$) on their way back to the limit cycle. This is not the case for the exact system, which is attracted to the \emph{nearest side} of the limit cycle.

\lfigure{FIG:XYHOPF_TRAJ01}{Sample trajectories computed for the Hopf normal form (\ref{EQ:LIN_APPROX_NOISY_LC_HOPF01}) (left column) and its linear approximation (\ref{EQ:LIN_APPROX_NOISY_LC_HOPF03}) (middle column), for weak noise ($\NSR=0.1$, top row) and strong noise ($\NSR=0.5$, bottom row). The right column shows the probability distribution of the $x$ component estimated for both processes. Note that in (c) both curves overlap to a great extent. In all cases $\alpha_0=\lambda=\alpha$. The kurtosis \citep{2005Kurtosis} in (f) is $\beta_2\approx 2.1$ for the exact process and $\beta_2\approx 2.6$ for its approximation.}{
\begin{tabular*}{\textwidth}{lll}
\includegraphics[width=0.32\textwidth,clip=true,trim=10mm 7mm 10mm 2mm]{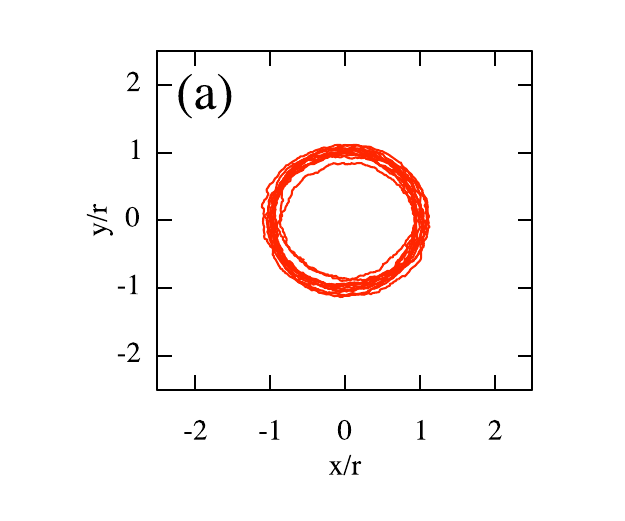}
&\includegraphics[width=0.32\textwidth,clip=true,trim=10mm 7mm 10mm 2mm]{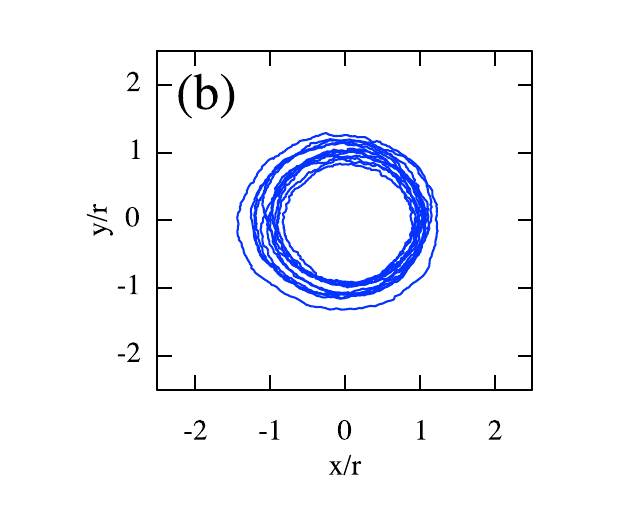}
&\includegraphics[width=0.32\textwidth,clip=true,trim=10mm 7mm 10mm 2mm]{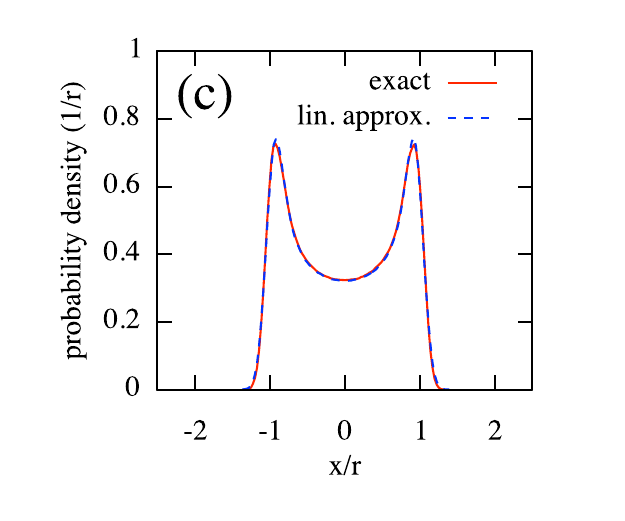}\\
\includegraphics[width=0.32\textwidth,clip=true,trim=10mm 7mm 10mm 2mm]{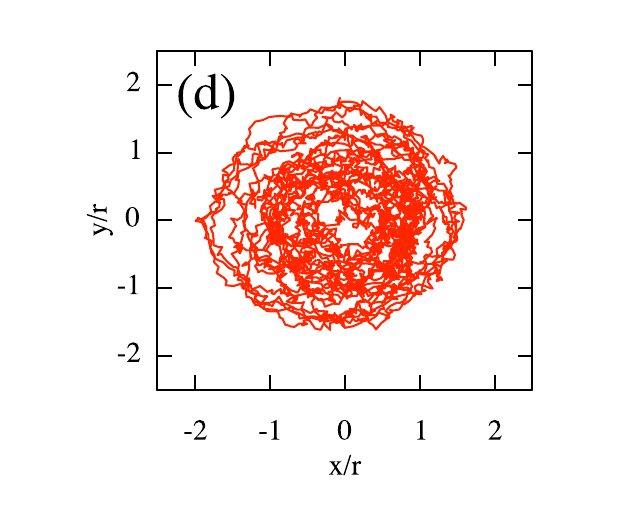}
&\includegraphics[width=0.32\textwidth,clip=true,trim=10mm 7mm 10mm 2mm]{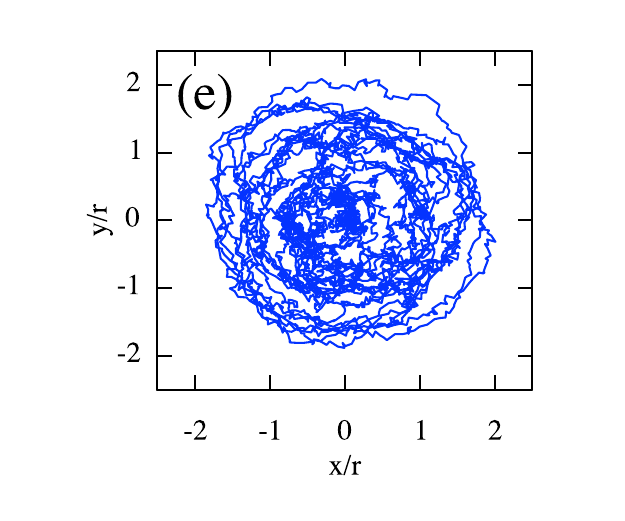}
&\includegraphics[width=0.32\textwidth,clip=true,trim=10mm 7mm 10mm 2mm]{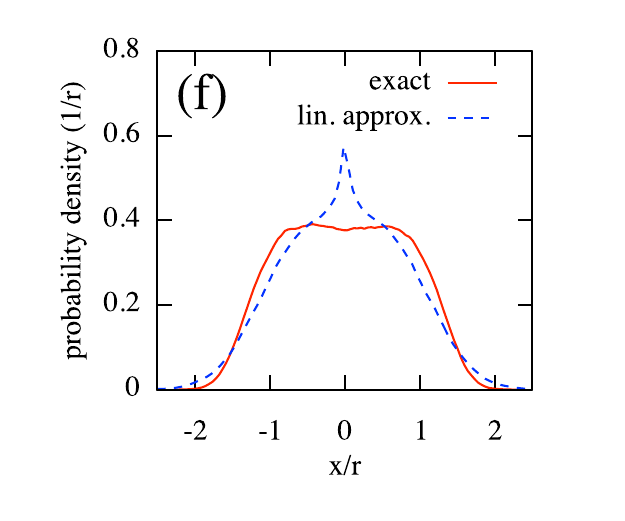}
\end{tabular*}
}

We compared the autocovariances ($\ACV$) and power spectral densities ($\PS$) for the first component $x$ of both processes, estimated from generated sample paths. We also compared the computed $\ACV$ and $\PS$ with the \emph{leading order formulas} \eqref{EQ:LEAD_ORD_APPR_ACV05} and \eqref{EQ:LC_SUPCRIT_HOPFBIF27}, respectively.
The linear approximation \eqref{EQ:LIN_APPROX_NOISY_LC_HOPF03}, as well as the leading order formulas, reproduce the exact $\ACV$ and $\PS$ to a great extent, whenever (i) the noise is sufficiently weak ($\NSR\lesssim 0.5$) and (ii) the focal and limit cycle frequencies are similar ($\abs{\alpha_0-\alpha}\lesssim\alpha/10$) (figures \ref{FIG:XYHOPF_PS_ACV01}a,b,d,e). 
The approximations fail when $\alpha$ differs significantly from $\alpha_0$ ($\abs{\alpha-\alpha_0}\gtrsim \alpha$) and $\NSR\gtrsim 0.1$ (figures \ref{FIG:XYHOPF_PS_ACV01}c,f). This is due to the nonlinear modulation of the phase speed at different deviations (rightmost term in \eqref{EQ:LIN_APPROX_NOISY_LC_HOPF01B}), not accurately represented in the linear approximation \eqref{EQ:LIN_APPROX_NOISY_LC_HOPF05}. In fact, the leading order formulas were derived by ignoring any deviation-dependent modulation of phase speed. Hence, one should expect a reduced accuracy of these approximations for systems with strong nonlinearities that modulate their phase speed far from the limit cycle.

\lfigure{FIG:XYHOPF_PS_ACV01}{Autocovariances ($\ACV$, plots (a,b,c)) and power spectral densities ($\PS$, plots (d,e,f)), computed for the Hopf normal form (\ref{EQ:LIN_APPROX_NOISY_LC_HOPF01}) (continuous line) and its linear approximation (\ref{EQ:LIN_APPROX_NOISY_LC_HOPF05}) (dashed line) using numerical simulations. The formulas (\ref{EQ:LEAD_ORD_APPR_ACV05}) for the $\ACV$ and (\ref{EQ:LC_SUPCRIT_HOPFBIF27}) for the $\PS$ are plotted for comparison (dotted line). Note that in figures (a,b,d,e) all three curves overlap to a great extent. Parameter values are $\alpha_0=\alpha$ and $\NSR=0.1$ for (a,d), $\alpha_0=\alpha$ and $\NSR=0.5$ for (c,e), $\alpha_0=\alpha/2$ and $\NSR=0.1$ for (c,f). In all cases $\lambda=\alpha$.}{
\begin{tabular*}{\textwidth}{ccc}
\includegraphics[width=0.325\textwidth,clip=true,trim=5mm 7mm 10mm 2mm]{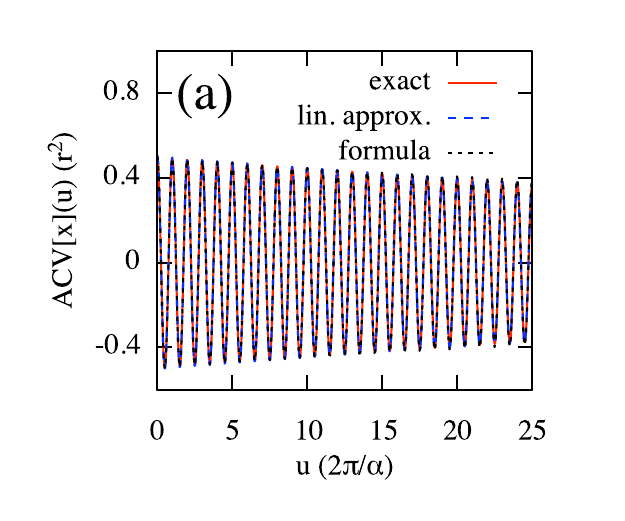}
&\includegraphics[width=0.325\textwidth,clip=true,trim=5mm 7mm 10mm 2mm]{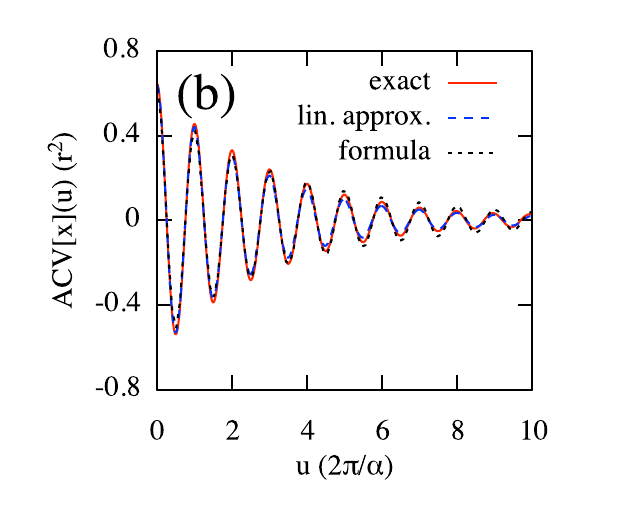}
&\includegraphics[width=0.325\textwidth,clip=true,trim=5mm 7mm 10mm 2mm]{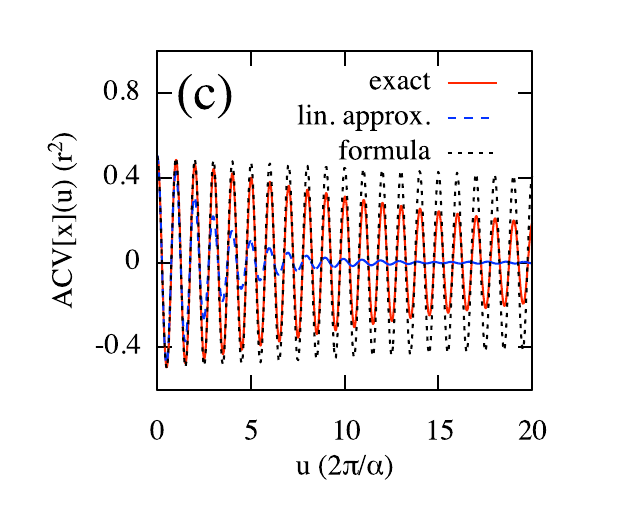}\\
\includegraphics[width=0.325\textwidth,clip=true,trim=5mm 7mm 10mm 2mm]{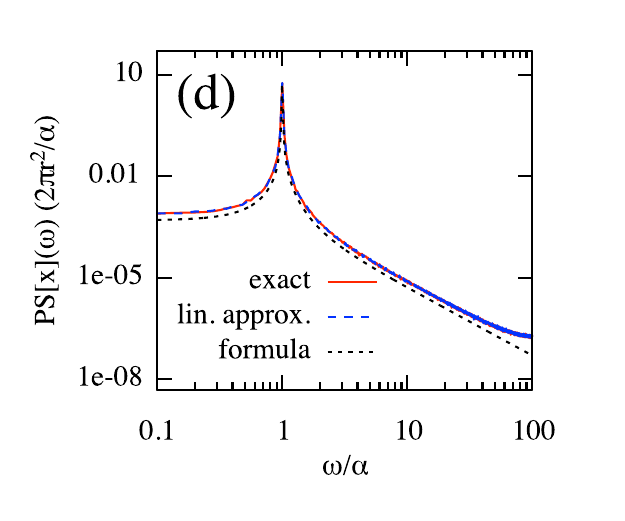}
&\includegraphics[width=0.325\textwidth,clip=true,trim=5mm 7mm 10mm 2mm]{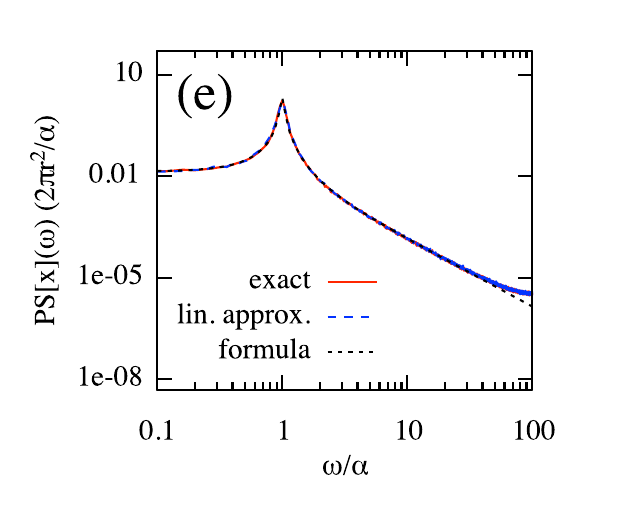}
&\includegraphics[width=0.325\textwidth,clip=true,trim=5mm 7mm 10mm 2mm]{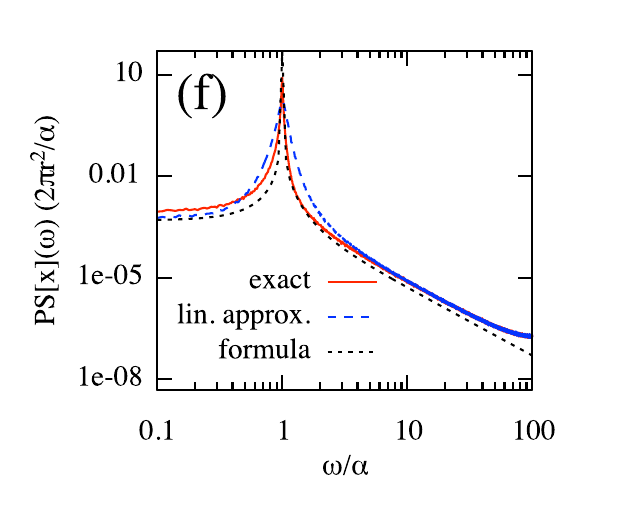}\\
\end{tabular*}
}

\section{Is the El-Ni\~no oscillation a noisy limit cycle?}
Decoherent periodic phenomena, such as animal population cycles, recurring disease outbreaks or climate oscillations, are conventionally described using noise-sustained oscillators \citep{McKane2005Predator-Prey,McKane2007Amplified,Greenwood2009Stochastic}, delay oscillators \citep{Battisti1989Interannual} or non-linear autoregressive models \citep{Ozaki1982The-statistical,Burgers1999The-El-Nino,Stephenson2000Is-the-North}. 
Given the ubiquity of deterministic ODE models exhibiting limit cycles, it is surprising that stochastically perturbed limit cycles have not seen an equally wide application in the interpretation of measured time series. In this article we have presented a generic way of understanding perturbed limit cycles using knowledge of the underlying deterministic dynamics, and have provided explicit formulas for the autocovariance and power spectrum that can be fitted to available time series.

As an illustrative example, we applied our results to the well studied but still poorly understood El Ni\~no Southern Oscillation (ENSO \citep{Trenberth1997The-Definition}). ENSO refers to quasi-periodic surface temperature variations in the tropical eastern Pacific Ocean recurring every 3--6 years, coupled to variations in air surface pressure in the tropical western Pacific.
The coupling between sea temperature and atmospheric pressure has been identified as a potential source of unstable modes \citep{Battisti1989Interannual,Cane1992Tropical}, however the decoherent nature of the ENSO has recently shifted the focus on models of oscillators sustained by stochastic weather processes \citep{Kleeman1994Limits,Penland1996A-stochastic,Moore1996The-dynamics,Moore1999Stochastic,Moore1999The-Nonnormal}.

Here we hypothesize that the ENSO can be described by an unstable mode giving rise to a limit cycle, which in turn is stochastically perturbed. Such a description reconciles past deterministic models with the observed irregularity of the ENSO and its lack of long-term predictability \citep{Kleeman1994Limits}. Independently fitting the ACV (Eq. \eqref{EQ:LEAD_ORD_APPR_ACV05}) and PS (Eq. \eqref{EQ:LC_SUPCRIT_HOPFBIF27}) to the ENSO time series (figures \ref{FIG:EL_NINO01}b,c) yields estimates for the deterministic Lyapunov exponent and noise variance. The fitted parameters are consistent between the ACV and PS, underlining the robustness of this approach. In particular, the ratio $\sigma^2/\ACV(0)$ is estimated at $0.83/\yr$ from the fitted ACV and at $0.96/\yr$ from the fitted PS, close to estimates obtained from noise-sustained oscillator models (approximately $1.41/\yr$) by \citet{Burgers1999The-El-Nino}. The fitted cycle period is about $4.2\ \yr$ in both cases. The focal Lyapunov exponent $\lambda/2$ is estimated at $0.15/\yr$ from the ACV and at $0.17/\yr$ from the PS, contrasting previous estimates using delay oscillator models ($1$--$1.5/\yr$) \citep{Battisti1989Interannual}.

\lfigure{FIG:EL_NINO01}{(a): Monthly anomaly of the El Ni\~no Southern Oscillation index N3.4 between January $1871$ and December $2007$ \citep{NCAR2013Nino}. (b,c): Autocovariance (b) and periodogram power (c), calculated using the time series data in (a). The dashed curves in (b) and (c) show the fitted autocovariance formula (\ref{EQ:LEAD_ORD_APPR_ACV05}) and fitted power spectrum formula (\ref{EQ:LC_SUPCRIT_HOPFBIF27}), respectively. Fitting was performed using least squares \citep{Bochkanov2013ALGLIB}.}{
\begin{tabular*}{\textwidth}{ccc}
\includegraphics[width=0.32\textwidth,clip=true,trim=13mm 6mm 10mm 10mm]{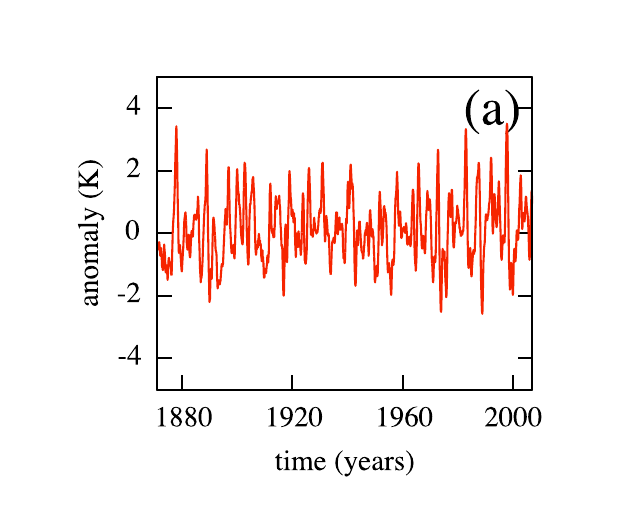}
& \includegraphics[width=0.32\textwidth,clip=true,trim=13mm 6mm 10mm 10mm]{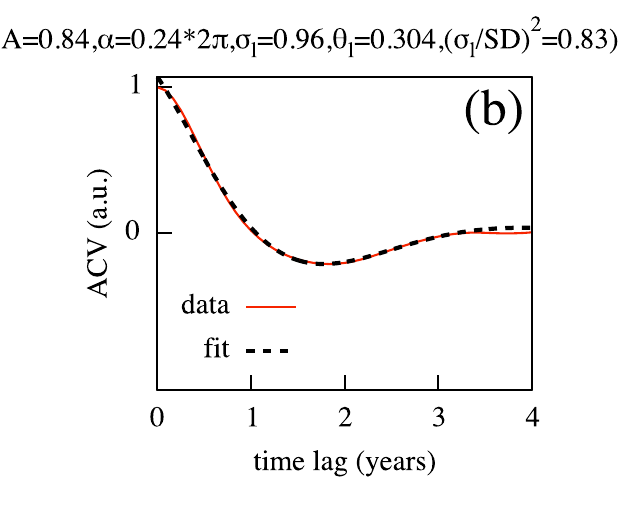}
& \includegraphics[width=0.32\textwidth,clip=true,trim=13mm 6mm 10mm 10mm]{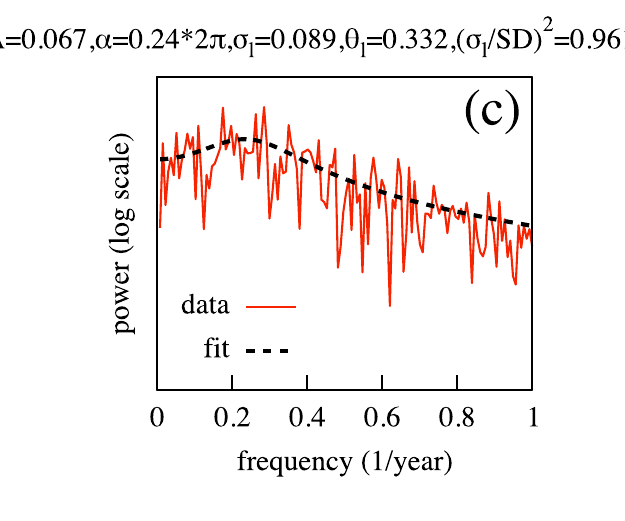}
\end{tabular*}
}

\section{Conclusions}
The approximation \eqref{EQ:LIN_APPROX_NOISY_LC15}, derived for the sample paths of noisy limit cycles, enables a qualitative understanding of perturbed systems exhibiting intrinsic cyclic dynamics. Our results suggest that such systems are better understood by separately considering the dynamics along and perpendicular to the deterministic limit cycle, as these are qualitatively different. While deviations from the limit cycle decay at a rate determined by the system's Jacobian, the phase drifts away from its deterministic value in a Brownian-motion-like manner. The rate of this random drift is, to leading order in the deviation, independent of the cycle's stability properties. This means that even systems with very stable limit cycles, can exhibit low temporal coherence. This becomes clear in the approximative formula \eqref{EQ:LEAD_ORD_APPR_ACV05} for the autocovariance of the Hopf normal form, which decays at an exponential rate that only depends on the cycle's radius and the noise strength. This decay makes noisy limit cycles fundamentally different from cyclostationary processes, which have a perfect time reference and whose autocovariance retains a non-decaying amplitude \citep{Jones1967Time,Gardner2006Cyclostationarity}. Therefore, noisy limit cycles are a potentially useful model for describing temporally decoherent cyclic processes such as animal population cycles \citep{Hornfeldt1994Delayed,Korpimaki2005Predatorinduced} or climate oscillations \citep{Trenberth1997The-Definition,Stephenson2000Is-the-North}, as exemplified above for the El Ni\~no Southern Oscillation. In fact, the formulas for the autocovariance and power spectrum allow, in principle, for a distinction of noisy limit cycles from noise-sustained oscillations. The latter also have well understood autocorrelations and power spectra \citep{McKane2005Predator-Prey,Tome2009Role}, and are often used to explain irregularly recurring phenomena \citep{Nisbet1976A-simple,Kaitala1996Cyclic,Moore1996The-dynamics,Moore1999Stochastic,Pascual2003Quasicycles,McKane2005Predator-Prey}.

The approximations derived in section \ref{SECTION:LC_HOPF}, both for the sample paths as well as the autocovariance and power spectrum, are strictly speaking only valid for systems with isotropic noise in the Hopf normal form \eqref{EQ:LIN_APPROX_NOISY_LC_HOPF01}. However, more general systems exhibiting stable limit cycles emerging through a Hopf bifurcation are expected to have autocovariances and power spectra that are qualitatively similar to the predictions given here, at least in the proximity of the bifurcation point. Hence, fitting these formulas to a given cyclic time series, as exemplified in this article, might provide a first estimate of the stability of the cycle and the amount of perturbations it is subject to. In contrast to conventional nonlinear autoregressive models \citep{Ozaki1982The-statistical}, the template formulas derived here permit a direct mechanistic interpretation of the estimated parameters.
In conclusion, this work provides a starting point for the qualitative understanding and statistical validation of stochastic differential equation models, exhibiting limit cycles in the deterministic limit.

\section{Acknowledgements}
This work was supported by the PIMS IGTC for Mathematical Biology, Canada. The authors would like to thank Priscilla Greenwood for comments.

\appendix

\section{Lemma on the existence of certain rotations}
\label{LEMMA:EX_CERT_ORTH_TRAFO}
Let $\bT(t)\in\R^n$ be a normalized vector that depends smoothly on time $t\geq0$. Let $P(t)$ be the hyperplane perpendicular to $\bT(t)$ and let $P_0=P(0)$, $\bT_0=\bT(0)$. Then there exists a unique family of orthogonal transformations $\U(t)\in \SO(n)$ ($t\geq0$), depending smoothly on time $t$, such that
\begin{enumerate}
\item \label{LEMMA:EX_CERT_ORTH_TRAFO:U00} $\U(0)$ is the identity,
\item \label{LEMMA:EX_CERT_ORTH_TRAFO:U01} $\U(t)P_0=P(t)$,
\item \label{LEMMA:EX_CERT_ORTH_TRAFO:U02} $\U(t)\bT_0=\bT(t)$,
\item \label{LEMMA:EX_CERT_ORTH_TRAFO:U03} $\frac{d \U}{dt}(t)\bz \perp P(t)$ for all $\bz\in P_0$.
\end{enumerate}
Moreover, $\U(t)$ satisfies the linear inhomogeneous differential equation
\lequation{EQ:EX_CERT_ORTH_TRAFO00}{
	\frac{d\U(t)}{dt} = -\bT(t)\frac{d\bT^\rT (t)}{dt}\U(t)\P_0 + \frac{d\bT (t)}{dt}\bT_0^\rT,
}
where $\P_0$ is the orthogonal projection onto $P_0$.

\paragraph{Proof} For notational simplicity we will denote $\bdot{X}=\frac{d}{dt}X$ for any time-dependent variable $X$. We start by showing the existence of $\U(t)$. Without loss of generality one can assume that
\lequation{EQ:EX_CERT_ORTH_TRAFO00B}{
	P_0=\curly{\bx=(x_1,..,x_n)\in\R^n : x_1=0}=\vspan\{\be_2,..,\be_n\}
}
and $\bT_0=\be_1$, where $\be_1,..,\be_n$ is the standard basis in $\R^n$.
Choose any $\A(t)\in\SO(n)$ depending smoothly on time and such that $\A(t)\bT_0=\bT(t)$ (such a transformation clearly exists). Denote $\bS(t)=\A^\rT(t)\bdot{\A}(t)\be_1$ and let 
\lequation{EQ:EX_CERT_ORTH_TRAF20}{
	\H(t) = \parmatrix{S_1(t) & -S_2(t) & \dots & -S_n(t)\\ S_2(t) & 0 & \dots & 0\\ \vdots & \vdots & \ddots & \vdots \\ S_n(t) & 0 & \dots & 0}, 
}
where $\bS=(S_1,..,S_n)^\rT$. Note that $S_1(t)=0$, since
\lequation{EQ:EX_CERT_ORTH_TRAFO22}{
	&\be_1^\rT \bS(t) = \be_1^\rT \A^\rT(t)\bdot{\A}(t)\be_1 = -\be_1^\rT\bdot{\A}^\rT(t)\A(t)\be_1 = - \block{\be_1^\rT \A^\rT(t)\bdot{\A}(t)\be_1}^\rT = -\be_1^\rT\bS(t).
}
In the 2nd step of \eqref{EQ:EX_CERT_ORTH_TRAFO22} we used the fact that $\A^\rT(t)\bdot{\A}(t) = -\bdot{\A}^\rT(t)\A(t)$, since $\A(t)$ is orthogonal. The matrix $\H(t)$ satisfies $\H^\rT(t) = -\H(t)$, $\H(t)\be_1 = \bS(t)$ and $\H(t)P_0\perp P_0$. Set $\D(t)=\A(t)\H(t)\A^\rT(t)$, then $\D(t)$ satisfies $\D^\rT(t) = -\D(t)$, $\D(t)\bT(t) = \bdot{\bT}(t)$ and $\D(t)P(t)\perp P(t)$.
Set $\U(t)$ as the solution to the ODE
\lequation{EQ:EX_CERT_ORTH_TRAFO24}{
	\bdot{\U}(t)=\D(t)\U(t), 
}
with initial value $\U(0)=\Id$. 
It is easy to see that $\bT(t)=\U(t)\bT_0$; indeed, $\U(t)\bT_0$ and $\bT(t)$ satisfy the same ODE and the same initial condition:
\lequation{EQ:EX_CERT_ORTH_TRAFO24A}{
	\frac{d}{dt}\parenth{\U\bT_0} = \bdot{\U}\bT_0=\D\cdot(\U\bT_0),\quad (\U\bT_0)(0) = \U(0)\bT_0 = \Id\bT_0 = \bT(0).
}
This proves claim (\ref{LEMMA:EX_CERT_ORTH_TRAFO:U02}). Moreover,
\lequation{EQ:EX_CERT_ORTH_TRAFO25}{
	\frac{d}{dt}\parenth{\U^\rT\U} = \bdot{\U}^\rT\U + \U^\rT\bdot{\U} = \U^\rT\D^\rT\U + \U^\rT\D\U
		= -\U^\rT\D\U + \U^\rT\D\U = 0,
}
that is, $\U^\rT(t)\U(t)=\U^\rT(0)\U(0)=\Id$ for all times $t$, implying $\U^\rT(t) = \U^{-1}(t)$. Hence, $\U(t)$ is indeed orthogonal. Therefore claim (\ref{LEMMA:EX_CERT_ORTH_TRAFO:U01}) follows from claim (\ref{LEMMA:EX_CERT_ORTH_TRAFO:U02}). Finally, for any $\bz\in P_0$ one has $\bdot{\U}(t) \bz = \D(t)\U(t)\bz\perp P(t)$, by property of $\D(t)$.

We shall now show \eqref{EQ:EX_CERT_ORTH_TRAFO00}, an immediate consequence of which will be the uniqueness of $\U(t)$. Denote $\V(t)=\bdot{\U}(t)$. Due to properties (\ref{LEMMA:EX_CERT_ORTH_TRAFO:U02}) and (\ref{LEMMA:EX_CERT_ORTH_TRAFO:U03}), $\V(t)$ can be written as
\lequation{EQ:EX_CERT_ORTH_TRAFO27}{
	\V(t) = \bT(t)\ba^\rT(t) + \bdot{\bT}(t)\bT_0^\rT
}
for some $\ba(t)\in P_0$. Since $\U(t)$ is orthogonal, we have $\U^\rT(t)\V(t)=-\V^\rT(t)\U(t)$. Hence, for any $\bz\in P_0$,
\lequation{EQ:EX_CERT_ORTH_TRAFO28}{
	(\ba^\rT(t)\bz)\bT_0=(\ba^\rT(t)\bz)\U^\rT(t)\bT(t) = \U^\rT(t)\V(t)\bz=-\V^\rT(t)\U(t)\bz = -\bT_0\bdot{\bT}^\rT(t)\U(t)\bz,
}
and thus $\ba=-\bdot{\bT}^\rT(t)\U(t)\P_0$. Therefore \eqref{EQ:EX_CERT_ORTH_TRAFO27} can be written as
\lequation{EQ:EX_CERT_ORTH_TRAFO29}{
	\V(t) = -\bT(t)\bdot{\bT}^\rT(t)\U(t)\P_0 + \bdot{\bT}(t)\bT_0^\rT,
}
which proves \eqref{EQ:EX_CERT_ORTH_TRAFO00}.\qed

\section{Lemma on the norm of certain rotations}
\label{APPENDIX:NORM_ROTATIONS}
Let $\bT(t)$, $P(t)$, $\bT_0$ and $P_0$ be as in \ref{LEMMA:EX_CERT_ORTH_TRAFO}. Let $\U(t)\in\SO(n)$ depend smoothly on time, such that $\U(t)\bT_0=\bT(t)$ and $\V(t)\bz\perp P(t)$ for all $\bz\in P_0$ (where $\V=\frac{d\U}{dt}$). Then $\norm{\V(t)}=\norm{\frac{d}{dt}\bT(t)}$.

\paragraph{Proof} For notational simplicity we will denote $\bdot{X}=\frac{d}{dt}X$ for any time-dependent variable $X$.
Since $\bT_0\perp P_0$ and $\bdot{\bT}(t)=\V(t)\bT_0$, we need to show that $\norm{\V(t)\bz}\leq \norm{\bz}\cdot\norm{\V(t) \bT_0}$ for all $\bz\in P_0$. Since $\V(t)\bz$ and $\bT(t)$ are parallel, one has
\lequation{EQ:NORM_ROTATIONS01}{
	\norm{\V(t)\bz} = \abs{\scalar{\V(t)\bz,\bT(t)}} = \abs{\scalar{\V(t)\bz,\U(t)\bT_0}} = \abs{\scalar{\U^\rT(t)\V(t)\bz,\bT_0}}.
}
Since $\U(t)$ is orthogonal, one has $\U^\rT(t)\V(t) = \V^\rT(t)\U(t)$. Hence
\lequation{EQ:NORM_ROTATIONS03}{
	\norm{\V(t)\bz} = \abs{\scalar{\U(t)\bz,\V(t)\bT_0}}\leq \norm{\U(t)\bz}\cdot\norm{\V(t)\bT_0} = \norm{\bz}\cdot \norm{\V(t)\bT_0},
}
as claimed.\qed

\section{Details on numerical methods}
\label{APPENDIX:DETAILS_NUM_METHODS}
For the generation of sample paths, we used an explicit two-step Runge-Kutta scheme of mean square order 3/2 \cite[\S 3.4, Theorem 3.3]{Milstein1995Numerical}, implemented in C++. 
 We normalized $r=1$ and $\alpha=2\pi$, and considered $\alpha_0\in[0.1,10]\cdot\alpha$ and $\lambda\in[0.1,10]\cdot\alpha$. Noise to signal ratios were considered in the range $\NSR\in[0.001,1]$. The integration time step was set to $10^{-4}$ cycle periods; decreasing it did not significantly change the outcome of the simulations. Autocovariances were estimated from the sample autocovariances of generated paths, spanning $10^6$ points over $10^4$ cycle periods \cite[\S 2.5.2]{Wei2005Time}. Power spectral densities were estimated by averaging the sample periodograms of $100$ independent sample paths, each spanning $10^4$ points over $10^2$ cycle periods \cite[\S 13.1]{Wei2005Time}. Probability distributions were estimated from $10^5$ points spanning $10^3$ cycle periods using a Gaussian kernel density estimator and Silverman's rule of thumb \citep[p. 48, eq. (3.31)]{Silverman1986Density}.



\end{document}